\documentclass[a4paper]{article}

\usepackage{arxiv}
\usepackage[utf8]{inputenc} 
\usepackage[T1]{fontenc}    
\usepackage{hyperref}       
\usepackage{url}            
\usepackage{booktabs}       
\usepackage{amsfonts}       
\usepackage{nicefrac}       
\usepackage{microtype}      
\usepackage{lipsum}		
\usepackage{graphicx}
\usepackage{natbib}
\usepackage{doi}

\usepackage{amsmath}
\usepackage{amsfonts}
\usepackage{amssymb}
\usepackage{amsthm}
\usepackage{newlfont}
\usepackage{amsbsy}
\usepackage{bm}
\usepackage[all]{xy}
\usepackage{textcomp}
\usepackage{color}
\usepackage{epsfig}
\usepackage{dsfont}
\usepackage{latexsym}
\usepackage{array}
\usepackage{longtable}
\usepackage{enumerate}
\usepackage{multicol}
\usepackage{mathrsfs}
\usepackage{cancel}
\usepackage{wasysym}

\newtheorem{teor}{Theorem}

\theoremstyle{definition}

\newtheorem*{ejem*}{Examples}
\newtheorem{defin}{Definition}

\theoremstyle{remark}

\newcommand{\complex}{\mathbf{\mathbb{C}}}
\newcommand{\reales}{\mathbf{\mathbb{R}}}

\newcommand{\natu}{\mathbf{\mathbb{N}}}
\newcommand{\llave}[1]{\left\{ #1\right\}}

\newcommand{\corch}[1]{\left[ #1\right]}
\newcommand{\paren}[1]{\left( #1\right)}

\title{Matrix bispectrality and noncommutative algebras: beyond the prolate spheroidals}
\date{\today}

\author{ 
\hspace{1mm}{F. Alberto Gr\"{u}nbaum}\\
Department of Mathematics\\
	University of California, Berkeley, CA\\
	\texttt{grunbaum@math.berkeley.edu} \\
	\And
\hspace{1mm}{Brian D. Vasquez} \\
	IMPA, Rio de Janeiro\\
	\texttt{bridava927@gmail.com} \\
		\And
\hspace{1mm}{Jorge P. Zubelli} \\
	Department of Mathematics, Khalifa University\\
	\texttt{zubelli@gmail.com} \\
}

\renewcommand{\underauthor}{Dedicated to Vaughan Jones, a truly inspirational friend, who left us some great memories.}
\renewcommand{\shorttitle}{}

\hypersetup{
pdftitle={A template for the arxiv style},
pdfsubject={q-bio.NC, q-bio.QM},
pdfauthor={David S.~Hippocampus, Elias D.~Striatum},
pdfkeywords={First keyword, Second keyword, More},
}

\begin{document}
\maketitle


\section{Preamble}

   The bispectral problem discussed in this paper has its origin in \cite{duistermaat1986differential}.
   It is motivated, as mentioned there, by an effort to understand and
   extend a remarkable phenomenon in Fourier analysis on the real line:
   the operator of time-and-band limiting is an integral operator admiting
   a second order differential operator with simple spectrum in its commutator.
   This property, which gives a good numerical way to compute the
   eigenfunctions of the integral operator,
   was put to good use in a series of papers by D. Slepian, H. Landau and 
   H. Pollak at Bells Labs back in the 1960's, see \cite{zbMATH03293058,zbMATH03293060,zbMATH03293059,zbMATH03293060,zbMATH03293061,zbMATH03588988,10.2307/2029386,Slepian1976OnB}
   and is of interest in one other contribution in this issue, see \cite{connes2021spectral,zbMATH05544468,connes2021prolate}.
    We are thankful to Luc Haine, of Louvain-la-Neuve, Belgium who alerted
one of us (AG) on Nov 2021 about a talk by Alain Connes in the series
Mathematicial Picture Language. This talk delivered on Dec 7th. 2021 can
be seen in Youtube and covers some of the contents of his work with various
collaborators on the zeta function and its relations to the prolate spheroidal
functions. For readers interested in this fascinating connection there is no
better way of learning about this than watching the lecture. See \cite{connes2021}.
 
   If one tries to extend this property beyond the Fourier case by adding
   a potential $V(x)$ to the operator $-\paren{\frac{d}{dx}}^2$ and replacing an expansion
   in exponentials by an expansion in term of the eigenfunctions $\psi(x,z)$
   of $-\paren{\frac{d}{dx}}^2+V(x)$ it appears plausible that a certain property of these
   eigenfunctions $\psi(x,z)$ could play a useful role. This property is now
   known as the bispectral property which will be formulated below. Its
   solution in the scalar case was the purpose of \cite{duistermaat1986differential}.

   The version of the bispectral problem that we discuss in this paper is a
   noncommutative one,
   obtained by allowing all objects in the original formulation to be
   matrix valued. The details are given later.

   Now for a bit of history: one of us (AG) gave a couple of
   talks, one at Vanderbilt
   in Nov 2013 and one in IMPA, Rio de Janeiro in March 2014 . The first talk was at
   the invitation
   of Vaughan Jones and started by saying that the topic was
   most likely of no interest
   to him. At some point in the talk Vaughan said "all of this
   is about bi-modules and subfactors". The second talk had both Vaughan
   and another one of us (JZ) in the audience, there was some more mention
   of these topics and then in \cite{Grunbaum2014SomeNM}, written
   in April 2014, a reference is made
   to a future joint paper with Vaughan
   "On the bimodule structure of the bispectral
   problem, in preparation". The occasion of this second talk was a visit that
   Vaughan (and for part of it with Wendy) did to Argentina, Uruguay and Brazil.
   In Buenos Aires he delivered one lecture, in Rio he delivered a series of
   talks but in Montevideo besides giving a lecture he was received in a
   private audience by President Jose Mujica, described as "the world humblest
   head of state" by Wikipedia. Vaughan whose command of the spanish language
   was quite high lamented that he had troubles understanding President Mujica.

   The examples in this
   short paper \cite{Grunbaum2014SomeNM} were discussed
   with Vaughan who showed some level of interest but
   we never managed to get him fully on board. Of course we assumed that
   we always had time to get him onto this project.
   The short paper \cite{Grunbaum2014SomeNM} makes three separate conjectures going with
   three different examples. These conjectures were proved to be correct in
   the thesis of the third of us (BVC) under the supervision of (JZ).
  
   And now that sadly we have no chance of benefiting from
   the insights that Vaughan
   would have brought in, we present the problem and some results to a wider
   audience in the hope that someone may surmise what it was that Vaughan had
   in mind. Getting someone involved in this effort would be a nice way of 
   honoring his memory.

\section{Commuting integral and differential operators}

 The bispectral problem, introduced in the next section is motivated by
the following very concrete problem in signal communication:
a signal of support in the interval $\corch{-T,T}$ is transmitted over a channel
that has bandwidth $\corch{-W,W}$ i.e., all frequencies in the signal beyond absolute value of $W$
cannot be sent over. A mathematical formulation is as follows: an arbitrary
signal in $L^2(\reales)$ is chopped to the interval $\corch{-T,T}$ and then its Fourier transform is chopped to the interval $\corch{-W,W}$. Denoting for simplicity these two
chopping operations by $T$ and $W$ we are dealing with the operator
$$E = W \mathcal{F} T$$
where $\mathcal{F}$ stands for the Fourier transform. The spectral analysis of this
operator, i.e., a look at its singular functions and singular values requires
the consideration of the operator $E^{*}E$. It is easy to see that $E^{*}E$ is an integral operator acting in 
$L^2 (-T,T)$ whose kernel is given by
$$\frac{\sin W (t-s)}{t-s} $$

and this bounded operator acts on a function in the space $L^2(-T,T)$ by 
\begin{equation*}
    (Kf)(s)=\int_{-T}^{T}\frac{\sin W (t-s)}{t-s} f(t)dt
\end{equation*}
for $f\in L^2(-T,T)$ and $s\in (-T,T)$. This $K$ commutes with the operator 
	$$(Lf)(x)=\paren{-\frac{d}{dx}}(T^2-x^2) \paren{\frac{df}{dx}} + W^2 x^2 f(x) $$
defined on $C^2$ functions.

One can show that this densely defined operator has a unique selfadjoint
extension in $L^2 (-T,T)$ with eigenfunctions and eigenvalues that depend on
the parameter $W$.

Its eigenfunctions are known as the prolate spheroidal wave functions, since this is one of the differential operators resulting in separating varables when solving for the eigenfunctions of the Laplacian on a prolate spheroid.

What other naturally appearing integral operators allow for commuting
differential operators?
Two other examples are the Bessel and Airy kernels, as in the work of Tracy and Widom, \cite{zbMATH00599692,zbMATH00509895} in the context of {\bf Random Matrix Theory.} For the Bessel case see also \cite{zbMATH03293061,zbMATH03783319}.
There are other examples, but the search is nowhere close to finished. 
The {\bf bispectral property}, to be formulated below was put
forward in \cite{duistermaat1986differential} as an important ingredient in the search
for more examples of this commuting property.

\section{The bispectral problem}

The problem was posed and solved in \cite{duistermaat1986differential}. It is as follows:

Find all nontrivial instances where a function $\psi(x,z)$ satisfies
\[
L\left( x,\frac {d}{dx} \right)\psi(x,z) \equiv \paren{-\paren{\frac{d}{dx}}^2 + V(x)}\psi(x,z)
= z\psi(x,z)
\]
as well as
\[
B\left(z,\frac {d}{dz} \right)\psi(x,z) \equiv \left( \sum_{i=0}^M
b_i(z)\left( \frac {d}{dz} \right)^i\right) \psi(x,z) =
\Theta(x)\psi(x,z).
\]
All the functions $V(x),b_i(k),\Theta(x)$ are, in principle, 
arbitrary except for
smoothness assumptions. Notice that here $M$ is arbitrary (finite). 
The operator $L$ could be of higher order, but in \cite{duistermaat1986differential} attention is restricted to order two.

\bigskip


\bigskip

The complete
solution is given as follows:

\bigskip
\noindent
{\bf Theorem } {\em If $M = 2$, then $V(x)$ is (except for
translation) either $\frac{c}{x^2}$ or $ax$, i.e.~we have a Bessel or an Airy case.
If $M > 2$, there are two families of solutions}

\begin{itemize}
\item[a)] {\em $L$ is obtained from $L_0 = -\paren{\frac{d}{dx}}^2$ by a finite number of
Darboux transformations $(L = AA^* \rightarrow {\tilde L} = A^*A)$.  In this
case $V$ is a rational solution of} the Korteweg-deVries hierarchy of equations. Here A is a first order differential operator.

\item[b)] {\em $L$ is obtained from $L_0 = -\paren{\frac{d}{dx}}^2 + \frac {1}{4x^2}$ after a
finite number of (rational) Darboux transformations.}
\end{itemize}

 In all cases we have a solution of the ad-conditions, a complicated system of nonlinear 
differential equations. These conditions are necesary and sufficient.

 Notice that the solutions organize themselves into nice manifolds.

\bigskip

The simplest example of case $a)$ follows from $L_0 = -\paren{\frac{d}{dx}}^2 $ by two Darboux transformations, one gets the operator
\[
L_2 = -\paren{\frac{d}{dx}}^2  + \frac {6(x^4+12 t_3 x)}{(x^3-t_3)^2}\ .
\]
In this case $\Theta(x) = x^4 - 4t_3 x$ and the differential operator in the
spectral parameter is

\[
B_2\paren{z,\frac{d}{dz}} = \left( -\paren{\frac{d}{dz}}^2 + \frac {6}{z^2} \right)^2 + 4 i t_3 \paren{\frac{d}{dz}}
\]

\bigskip

The potential in the operator $L_2= -\paren{\frac{d}{dx}}^2+ V(x,t)$ above satisfies the KdV equation.

It was later observed by Magri and Zubelli, see \cite{zbMATH00028878} that in case $b)$ we are dealing
with rational solutions of the Virasoro equations (i.e. master symmetries of KdV). 
The bigger picture became more apparent in the work \cite{zbMATH01514630} where it is shown that the generic rational potentials that decay at infinity and remain rational by all the flows of the master-symmetry KdV hierarchy are bispectral potentials for the Schrödinger operator.

In case $a)$ the space of common solutions has dimension one and in case $b)$ it has dimension two.  One refers to these as the rank one and rank two situations.

Observe that the ``trivial cases'' when $M=2$ are self-dual in the sense that
since the eigenfunctions $\psi(x,z)$ are functions either of the product $xz$ or
of the sum $x+z$, one gets $B$ by replacing $z$ for $x$ in $L$.
The {\it bispectral involution} introduced in \cite{zbMATH00438440} shows how this
can be adapted in the "higher order cases".

\section{The non-commutative version of the bispectral problem}

A first noncommutative (or matrix) version of the bispectral problem was
considered in J. Zubelli Ph.D. thesis at Berkeley, see also \cite{zubelli1990differential,Zubelli1992a}\cite{zbMATH00040257} in the situation where both the physical space and spectral operators act on the same side of the eigenfunction and the eigenvalues are both
scalar valued. Later on, several other versions were considered. See \cite{zbMATH01694217,zbMATH06530535,zbMATH05346328,zbMATH02009271,Grunbaum2014SomeNM,zbMATH06722531} and references therein. 
The noncommutative version of the bispectral problem displayed interesting connections with soliton equations as well. Indeed, in \cite{zbMATH00040257} it was shown that a large class of rational solutions to the AKNS hierarchy \cite{zbMATH03634788} led to matrix differential operators that displayed the bispectral property. Among the important equations in Mathematical Physics that are covered by the AKNS hierarchy one finds the modified KdV and the nonlinear Schr\"{o}dinger equation. 
The matrix differential operator that appeared in this case was in turn related to Dirac operators.
The connection between bispectrality and another important topic in Mathematical Physics, namely Huygens' principle in the strict sense \cite{BerestVeselov94} turned out to appear also in the context of Matrix Bispectrality. Indeed, in \cite{MR2201201, MR2177466, MR2091493, Chalub2001}, it was shown 
that rational solutions to the AKNS hierarchy led to Dirac operators which satisfy Huygens' principle in the stric sense. In other words, the fundamental solutions of the perturbed Dirac equation in a suitably high space-time dimension had its support precisely on the surface of the light cone and not in its interior. 
Another interesting connection between Matrix bispectrality to soliton equations of Mathematical Physics was explored in \cite{zbMATH01694217}.

In this present paper 
we take the {\bf bi-module} structure of the problem into account and let the operators act on different sides as well as allow both eigenvalues to be matrix valued.

\bigskip

We consider triplets $(L,\psi, B)$ satisfying the equations 
	\begin{equation}\label{bispec}
	L\psi(x,z)=\psi (x,z)F(z) \hspace{1 cm}  (\psi B) (x,z)=\theta(x)\psi(x,z)
	\end{equation}
	with $L=L(x,\frac{d}{dx})$, $B=B(z,\frac{d}{dz})$  linear matrix  differential operators, i.e., $L\psi=\sum_{i=0}^{l}a_{i}(x) \cdot \paren{\frac{d}{dx}}^i \psi$,
	$\psi B=\sum_{j=0}^{m}\paren{\frac{d}{dz}}^j \psi \cdot b_{j}(z)$. The functions $a_{i}:U \subset \complex \rightarrow M_{N}(\complex)$, $b_{j}:V\subset \complex \rightarrow M_{N}(\complex),$ $F:V\subset \complex \rightarrow M_{N}(\complex),$ $\theta:U\subset \complex\rightarrow M_{N}(\complex)$ and the nontrivial common eigenfunction 
	$\psi:U\times V\subset \complex^2\rightarrow M_{N}(\complex)$ are in principle compatible sized meromorphic matrix valued functions defined in suitable open subsets $U,V\subset \complex$. 
	%


A triplet $(L,\psi, B)$ satisfying \eqref{bispec} is called a bispectral triplet. 

The study of the structure of the algebra of possible $\theta(x)$
going with a fixed bispectral $\psi(x,z)$ was first raised
in \cite{zbMATH05249007} and analyzed in \cite{zbMATH05896191,zbMATH05199815}.
See also \cite{zbMATH06859910} and \cite{zbMATH07004359}.

We consider now the examples and conjectures given in \cite{Grunbaum2014SomeNM} as well as
their validation and further description given in \cite{VZ1}. 

For the benefit of the reader we give a few definitions before giving some explicit results in the next section.

\begin{defin}
Let  $\mathbb{K}$ be a field, $C$ be a $\mathbb{K}$-algebra, $A$ a subring of $C$ and $S\subset C$. We define 
$$A\cdot <S>=span_{\mathbb{K}}\llave{\prod_{j=1}^{n}s_j\mid s_1, ... , s_n \in S\cup A, n \in \natu } \mbox{ ,}$$
where the noncommutative product is understood from left to right, i.e., 
$ \prod_{j=1}^{n+1}s_j := (\prod_{j=1}^{n}s_j )s_{n+1}, $ for $n=0,1,2,\cdots.$ For completion,
$\prod_{j=1}^{0}s_j:=1.$
\end{defin} 

The set $A\cdot<S>$ is called the {\em subalgebra generated by} $S$ over $A$ and we call an element $f\in A\cdot<S>$ a noncommutative polynomial with 
coefficients in $A$ and set of variables $S$. 

\begin{defin}
    Let $C$ be a noncommutative ring and $A$ a subring of $C$. We say that an element $\alpha \in C$ is integral over $A$ if there exists a 
    noncommutative polynomial $f$ with coefficients in $A$ such that $f(\alpha)=0$. Furthermore, we say that $\beta\in C $ is integral 
    over $\alpha \in C$ if $\beta$ is integral over $A\cdot<\alpha>$. Finally, $\alpha$ and $\beta$ are associated integral if  $\alpha$ is integral over $\beta$ and $\beta$ is integral over $\alpha$.
\end{defin}
In order to characterize the algebraic structure of bispectrality in 
the present noncommutative context, we start with the following definitions.

\begin{defin}
Let  $\mathbb{K}$ be a field, we denote by $\mathbb{K}\langle x_{\lambda}\mid \lambda\in \Lambda  \rangle$ the free algebra generated by the letters $x_{\lambda}$, $\lambda\in \Lambda$ i.e., $$\mathbb{K}\langle x_{\lambda}\mid \lambda\in \Lambda  \rangle=
\bigoplus_{F\subset \Lambda, F \text{finite}} \bigoplus_{\lambda\in F}\mathbb{K}\cdot x_{\lambda}.$$
\end{defin}
\begin{defin}
        Let $A$ be a $\mathbb{K}$-algebra. A presentation for an algebra $A$ is a triple $(\mathbb{K}\langle x_{\lambda}\mid \lambda\in \Lambda  \rangle,f,I)$ such that $I\subset A$ is an ideal and $f:\mathbb{K}\langle x_{\lambda}\mid \lambda\in \Lambda  \rangle/ I \rightarrow A$ is an isomorphism. Furthermore, we say that $A$ is finitely generated if there exists a presentation with $\Lambda$ finite and finitely presented if there exists a presentation with $\Lambda$ finite and the ideal $I$ is generated by finitely many elements.
\end{defin}

\section{ Three examples }

Take for $\Psi(x,z)$ the matrix valued function
\[
\Psi(x,z) = e^{xz}\begin{pmatrix}
z - 1/x & 1/x^2 \\
0 & z - 1/x
\end{pmatrix}
\]
and consider all instances of matrix-valued polynomials $\theta(x)$ and differential operators $B$ (with matrix coefficients $b_i(z)$) such that
\[
\Psi B \equiv \sum_{i=1}^m \paren{\frac{d}{dz}}^i\Psi\ b_i = \theta(x)\Psi(x,z).
\]

In this case one has
\[
L\Psi = -z^2\Psi
\]
with
\[
L = -\paren{\frac{d}{dx}}^2 + 2\begin{pmatrix}
1/x^2 & -2/x^3 \\
0 & 1/x^2
\end{pmatrix}.
\]

In other words for this specific differential operator in the variable $x$ we are asking for all bispectral "partners" of $L$.

\bigskip

One finds that one such pair $(B,\theta)$ is given by
\[
B = \paren{\frac{d}{dz}}^3 - 3\paren{\frac{d}{dz}}^2 \cdot \frac {1}{z} + 3\paren{\frac{d}{dz}}\cdot \frac {1}{z^2} + 3\begin{pmatrix}
0 & 1/z^2 \\
0 & 0
\end{pmatrix}
\]
and $\theta(x)$ the scalar-valued polynomial
\[
\theta(x) = x^3.
\]

The set of all possible $\theta(x)$ is given by the following
subalgebra $\mathbb{A}$. The complete statement is given by 

\bigskip



\begin{teor}\label{first1}
	Let $\Gamma$ be the sub-algebra of $M_{2}(\complex)[x]$ of the form 
	\begin{equation*}
	\paren{
		\begin{matrix} 
		r_{0}^{11} & r_{0}^{12} \\
		0 & r_{0}^{11} 
		\end{matrix}}+  \paren{
		\begin{matrix} 
		r_{1}^{11} & r_{1}^{12} \\
		0 & r_{1}^{11} 
		\end{matrix}}x+ \paren{
		\begin{matrix} 
		r_{2}^{11} & r_{2}^{12} \\
		r_{1}^{11} & r_{2}^{22} 
		\end{matrix}}x^{2}+
	\paren{
		\begin{matrix} 
		r_{3}^{11} & r_{3}^{12} \\
		r_{2}^{22}+r_{2}^{11}-r_{1}^{12} & r_{3}^{22} 
		\end{matrix}}x^3+x^{4}p(x),
	\end{equation*} 
	where $p\in M_{2}(\complex)[x]$ and all the variables $r_{0}^{11},r_{0}^{12}, r_{1}^{11},r_{1}^{12},r_{2}^{11},r_{2}^{22},r_{3}^{11},r_{3}^{12}, r_{3}^{22} \in \complex$. Then $\Gamma=\mathbb{A}$. Moreover, for each $\theta$ we have an explicit expression for the operator $B$.
	
	Furthermore, we have the presentation $\mathbb{A}=\complex\cdot \langle \alpha_0, \alpha_1 \mid I=0\rangle$ with
the ideal $I$ given by 
\begin{equation*}
I :=\langle \alpha_{0}^2,\alpha_{1}^3+\alpha_{0}\alpha_{1}\alpha_{0}-3\alpha_{1}\alpha_{0}\alpha_{1}+\alpha_{0}\alpha_{1}^2+\alpha_{1}^2\alpha_{0} \rangle \mbox{ .}
\end{equation*}
\end{teor}
This is an example of an algebra with an integral element over a nilpotent one.

For the next example take for $\psi(x,z)$ the matrix-valued function
\[
\psi(x,z) = \left[\frac{d}{dx} - \begin{pmatrix}
1/x & -1/x^2 & 1/x^3 \\
0 & 1/x & -1/x^2 \\
0 & 0 & 1/x
\end{pmatrix} \right]\ e^{xz}I = e^{xz} \begin{pmatrix}
z - 1/x & 1/x^2 & -1/x^3 \\
0 & z - 1/x & 1/x^2 \\
0 & 0 & z - 1/x
\end{pmatrix}.
\]
Here one can see that
\[
L\psi = -z^2\psi
\]
with
\[
L = -\paren{\frac{d}{dx}}^2 + 2\begin{pmatrix}
1/x^2 & -2/x^3 & 3/x^4 \\
0 & 1/x^2 & -2/x^3 \\
0 & 0 & 1/x^2
\end{pmatrix}.
\]

The results in this case about the set of all possible $\theta(x)$ are given 
below.

\begin{teor}\label{second}

	Let $\Gamma$ the sub-algebra of $M_{3}(\complex)[x]$ of the form 
	\begin{equation*}
	\paren{
		\begin{matrix} 
		r_{0}^{11} & r_{0}^{12} & r_{0}^{13} \\
		0 & r_{0}^{22} & r_{0}^{23} \\
		0 & 0 & r_{0}^{11}
		\end{matrix}}+ \paren{
		\begin{matrix} 
		r_{1}^{11} & r_{1}^{12} & r_{1}^{13} \\
		r_{0}^{22}-r_{0}^{11} & r_{1}^{22} & r_{1}^{23} \\
		0 & r_{0}^{22}-r_{0}^{11} & r_{1}^{11}+r_{0}^{23}-r_{0}^{12}
		\end{matrix}}x
	\end{equation*}
	\begin{equation*}
	+ \paren{
		\begin{matrix} 
		r_{2}^{11} & r_{2}^{12} & r_{2}^{13} \\
		r_{1}^{22}-r_{1}^{11}-r_{0}^{23}+r_{0}^{12} & r_{2}^{22} & r_{2}^{23} \\
		r_{0}^{22}-r_{0}^{11} & r_{1}^{22}-r_{1}^{11} & r_{2}^{11}+r_{1}^{23}-r_{1}^{12}
		\end{matrix}}x^{2}
	+  \paren{
		\begin{matrix} 
		r_{3}^{11} & r_{3}^{12} & r_{3}^{13} \\
		r_{3}^{21}& r_{3}^{22} & r_{3}^{23} \\
		r_{1}^{22}-2r_{1}^{11}-r_{0}^{23}+r_{0}^{12} & r_{3}^{32} & r_{3}^{33}
		\end{matrix}}x^{3}   
	\end{equation*}
	\begin{equation*}
	+ \paren{
		\begin{matrix} 
		r_{4}^{11} & r_{4}^{12} & r_{4}^{13} \\
		r_{4}^{21} & r_{4}^{22} & r_{4}^{23} \\
		r_{3}^{32}+r_{3}^{21}-r_{2}^{22}-r_{2}^{11}+r_{1}^{12} & r_{4}^{22} & r_{4}^{33}
		\end{matrix}}x^{4}
	\end{equation*}	
	\begin{equation*}
+  \paren{
	\begin{matrix} 
	r_{5}^{11} & r_{5}^{12} & r_{5}^{13} \\
	r_{5}^{21}& r_{5}^{22} & r_{5}^{23} \\
	r_{4}^{32}+r_{4}^{21}-r_{3}^{33}-r_{3}^{22}-r_{3}^{11}+r_{2}^{23}+r_{2}^{12}-r_{1}^{13} & r_{5}^{32} & r_{5}^{33}
	\end{matrix}}x^{5}	+x^{6}p(x) \mbox{ , }
	\end{equation*} 
	where $p\in M_{3}(\complex)[x]$ and all the variables $r_{0}^{11}, r_{0}^{12},...,r_{5}^{33}\in \complex$ are arbitrary.
	
	Then, $\Gamma =\mathbb{A}$ and for each $\theta$ we have  an explicit expression for the operator $B$.
	
	Furthermore, we have the presentation $\mathbb{A}=\complex\cdot \langle \alpha_2, \alpha_3 \mid I=0\rangle$ with
	\begin{equation*}
	I=\langle\alpha_{2}^3,\alpha_{3}^2-\alpha_{3}, (\alpha_{3}\alpha_{2})^2 \alpha_{3}-4\alpha_{3}\alpha_{2}^2 \alpha_{3}\rangle \mbox{ .}
	\end{equation*}
\end{teor}
This is an example of an algebra with nilpotent and idempotent associated elements.


As the last example we consider a case when both "eigenvalues" $F$ and $\theta$ are matrix valued. 
Let 
\begin{equation*}
\psi(x,z)=\frac{e^{xz}}{(x-2)xz}\paren{\begin{matrix} 
\frac{x^3z^2-2x^2z^2-2x^2z+3xz+2x-2}{xz} & \frac{1}{x} \\
\frac{xz-2}{z} & x^2z-2xz-x+1
\end{matrix}}
\end{equation*}
and \begin{equation*}
L=\paren{\begin{matrix} 
     0 & 0 \\
	0 & 1
	\end{matrix}}.\paren{\frac{d}{dx}}^2
+\paren{\begin{matrix} 
0 & \frac{1}{(x-2)x^2} \\
	-\frac{1}{x-2} & 0
	\end{matrix}}.\paren{\frac{d}{dx}}
+ \paren{\begin{matrix} 
	-\frac{1}{x^2(x-2)^2} & \frac{x-1}{x^3(x-2)^2} \\
	\frac{2x-1}{x(x-2)^2} & -\frac{2x^2-4x+3}{x^2(x-2)^2}
	\end{matrix}},
\end{equation*}
then $L\psi=\psi F$ with 
\begin{equation*}
F(z)=\paren{\begin{matrix} 
	0 & 0 \\
	0 & z^2
	\end{matrix}}.
\end{equation*}
It  is easy to check that $\psi B=\theta \psi$ for 
\begin{equation*}
B=\paren{\frac{d}{dz}}^3.\paren{\begin{matrix} 
                   0 & 0 \\
	               1 & 0
	\end{matrix}}+ \paren{\frac{d}{dz}}^2.
\paren{\begin{matrix} 
	0 & 0 \\
	-\frac{2z+1}{z} & 0
	\end{matrix}}+
	\paren{\frac{d}{dz}}.\paren{\begin{matrix} 
		1 & 0 \\
		\frac{2(z-1)}{z^2} & 1
		\end{matrix}}
		+\paren{\begin{matrix} 
			-z^{-1} & 0 \\
			6z^{-3} & z^{-1}
			\end{matrix}}
\end{equation*}
and 
\begin{equation*}
\theta(x)=
\paren{\begin{matrix} 
	x & 0 \\
	x^2(x-2) & x
	\end{matrix}}.
\end{equation*}
In this case we characterize the algebra $\mathbb{A}$ of all polynomial $F$ such that there exist $L=L\paren{x,\frac{d}{dx}}$ with $L\psi=\psi F$ as follows

\begin{teor}\label{calogero}
	Let $\Gamma$ be the sub-algebra of $M_{2}(\complex)[z]$ of the form 
	\begin{equation*}
	\paren{\begin{matrix} 
		a & 0 \\
		b-a & b
		\end{matrix}}+
	\paren{\begin{matrix} 
		c & c \\
		a-b-c & -c
		\end{matrix}}z+
	\paren{\begin{matrix} 
		a-b-c & c+a-b\\
		d & e
		\end{matrix}}\frac{z^2}{2}+z^3p(z),
	\end{equation*}
	where $p\in M_{2}(\complex)[z]$ and all the variables $a,b,c,d,e$ are arbitrary. Then $\Gamma=\mathbb{A}$.
	
		Furthermore, we have the presentation $\mathbb{A}=\complex\cdot \langle \theta_{1}, \theta_{3},\theta_{4}, \theta_{5}\mid I=0\rangle$ with
	
	\begin{equation*}
	I=\langle\theta_{1}^2-\theta_{1},\theta_{4}^2,\theta_{4}\theta_{5},
	\theta_{4}\theta_{1}+\theta_{4}\theta_{3}-2\theta_{4}-\theta_{5}\theta_{4}-\theta_{5}^2,
	\theta_{3}^2-\theta_{3}+\theta_{5}-3\theta_{3}\theta_{4}\theta_{3}\theta_{5}-\theta_{1}\theta_{4}-\theta_{5}\theta_{1},
		\end{equation*}
	\begin{equation*}
	\theta_{3}\theta_{1}-\theta_{1}-\theta_{4}-\frac{1}{2}\theta_{4}\theta_{1}+\frac{1}{2}\theta_{4}\theta_{3}+\theta_{5}\theta_{1}-\frac{1}{2}\theta_{5}\theta_{4}
	+\frac{1}{2}\theta_{5}^2+\theta_{3}\theta_{4}-\theta_{1}\theta_{5}-\theta_{3}\theta_{5},
		\end{equation*}
	\begin{equation*}
	\theta_{1}\theta_{3}-\theta_{3}+\theta_{4}+\theta_{5}-\frac{3}{2}\theta_{4}\theta_{1}+\frac{3}{2}\theta_{4}\theta_{3}-2\theta_{5}\theta_{1}-\frac{3}{2}\theta_{5}\theta_{4}
	+\frac{3}{2}\theta_{5}^2+3\theta_{3}\theta_{4}+\theta_{3}\theta_{5},
		\end{equation*}
	\begin{equation*}
	\theta_{5}\theta_{3}-\theta_{4}\theta_{1}+\theta_{4}\theta_{3}-\theta_{5}\theta_{1}-\theta_{5}\theta_{4}+\theta_{5}^2,
	\theta_{5}\theta_{1}\theta_{5}-\theta_{5}^2\theta_{1}-\theta_{5}\theta_{4},
	\theta_{5}\theta_{4}\theta_{1}-\theta_{5}^3+\theta_{5}\theta_{1}\theta_{4}+\theta_{5}^2\theta_{1},
		\end{equation*}
	\begin{equation*}
	\theta_{4}\theta_{1}\theta_{5}+\theta_{4}\theta_{3}\theta_{5}-\theta_{3}^3,
	\theta_{5}\theta_{3}\theta_{4}+\theta_{5}\theta_{1}\theta_{4}
	\rangle.
	\end{equation*}
	
\end{teor}

This is an example of an algebra with two integer elements over one nilpotent and one idempotent. This is linked to the spin Calogero  systems whose relation with bispectrality can be found in \cite{zbMATH05577256,zbMATH06530535}.

Theorems \ref{first1}, \ref{second} and \ref{calogero} give positive answers to the Conjectures 1, 2 and 3  of \cite{Grunbaum2014SomeNM} about three bispectral full rank 1 algebras. Moreover, these algebras are Noetherian and finitely generated because they are contained in the $N\times N$ matrix polynomial ring $M_{N}(\mathbb{K}[x])$ (See \cite{VZ1}).

We close this section by remarking the important role played by Darboux transformations,
which goes back to \cite{duistermaat1986differential} and \cite{Zubelli1991}. Indeed, 
the three examples in this section are instances of rational Darboux
transformations from the scalar matrix exponential functions. All such
Darboux transformations were shown to be bispectral in Theorems~1.1
and 1.2 of  \cite{zbMATH06722531}.

\section{A few extensions of the problems discussed above}
Since one of the goals of this paper is to serve as an ``invitation'' to
look at this problem extended to a wider audience, we give a road map with
some selected references.
Indeed, the bispectral problem has many different incarnations, and in our opinion we are still far from having a unified theory. 

In the scalar case here are some early papers that should be mentioned are \cite{zbMATH01000896,zbMATH00438440,zbMATH01186161,zbMATH01410877,zbMATH01213907,zbMATH01425111,zbMATH01149951,zbMATH01818998,zbMATH06722531}.

One natural issue concerns the numerical aspects involving the  prolate spheroidal functions. In this case, the reader may want to consult the book \cite{zbMATH06189114}.

Another direction concerns, the purely discrete (actually finite) version of time-band-limiting. 
In an effort to better understand the commuting property in question, one of us
looked at the case when the real line is replaced by the N roots of unity.
See \cite{zbMATH03822507} as well as \cite{zbMATH03877750,zbMATH03998232}.
 Expression (11) in \cite{zbMATH03822507} contains a small typo: the $r(2)$ in the denominator should
be replaced by $r(1)$.

Moving on to the discrete-continuous version of the bispectral problem, we have that for the scalar case, involving orthogonal polynomials satisfying differential equations the problem had already been considered by S. Bochner (and even earlier).
A very good introduction to this is given in \cite{zbMATH01684158} and its references. See also \cite{zbMATH01684157, zbMATH01513050,zbMATH01121135,MR2000k:33017}.

For the matrix valued case there are two sources of early examples, one resulting from the theory of matrix valued spherical functions,
see \cite{zbMATH01709246,zbMATH01763777,zbMATH01911954,zbMATH02077009,zbMATH02177602,zbMATH02230067}, and another
one see \cite{zbMATH02207651,zbMATH05202123,zbMATH05351608}. See also \cite{zbMATH02009271} as well as \cite{zbMATH05346328,zbMATH02210918,1988CMaPh.119..385R}.

Solutions of the bispectral  problem can be used to obtain integral operators which reflect some ordinary differential operator in the sense of (for instance) \cite{zbMATH07315538}. This fact generalizes the commuting property in the scalar case. It would be interesting to see whether this could be extended to the matrix case.

In the (noncommutative) matrix case, by considering operators in the physical and spectral variables acting from opposite directions, we maintain the Ad-conditions that played a substantial role in \cite{duistermaat1986differential}. In this case, this leads to the embedding of the bispectral algebras of eigenvalues into the matrix polynomial algebra $M_{N}(\complex[x])$. See \cite{VZ1}.
    
Another natural direction would be to look for a characterization for algebras relating to the spin-Calogero system for matrices of arbitrary size of matrix $N$. The examples $N=1$ and $N=2$ were generalized to arbitrary matrix size $N$ and was characterized as a subalgebra of $M_{N}(\complex[x])$ using a family of maps $\llave{P_{k}}_{k\in \natu}$ satisfying some nice properties such as translation and a product rule similar to the Leibniz rule. See \cite{VZ1}.


Finally, let us go back to the original problem.
We recall that both in the scalar and in the matrix case the motivation behind
the bispectral problem was a desire to understand "what is behind" the remarkable commutativity property between the operator of time-and-band limiting in the
Fourier, Bessel and Airy cases.
It was hard to suspect that this problem would have connections with many of the recent developments in integrable systems. There has been some progress in connecting the bispectral problem with the commutativity property mentioned above, and here again-at least in the scalar case- there are some connections with Integrable systems, see for instance \cite{grun2021,zbMATH07244600,zbMATH07315538,casper2020reflective,CGYZ3}. For connections between the bispectral problem and the commuting property see \cite{zbMATH06813265,zbMATH06691656,zbMATH06979371}.

\section*{Acknowledgments}

BVDC was supported by CAPES grants 88882 332418/2019-01 as well as IMPA.
JPZ was supported by CNPq grants 302161 and 47408, as well as by FAPERJ under the program \textit{Cientistas do Nosso Estado} grant E-26/202.927/2017, Brazil. BDVC and JPZ acknowledge the support from the FSU-2020-09 grant from Khalifa University, UAE.

\end{document}